\documentclass{article}
\usepackage{amsmath}
\usepackage{amsfonts}
\usepackage{amsthm}
\begin{document} 

\def\l2e{\langle\!\langle}
\def\r2e{\rangle\!\rangle}
\def\la{\lambda}
\def\Bl2e{\Bigl\langle\!\Bigl\langle}
\def\Br2e{\Bigr\rangle\!\Bigr\rangle}
\newtheorem{theorem}{Theorem}
\newtheorem{comment}{Remark}
\newtheorem{corollary}{Corollary}
\newtheorem{propo}{Proposition}

\begin{titlepage} 
\title{Completion of a Rational Function Sequence of Carlitz}
\author{ Leonard M. Smiley\\
Department of Mathematics\\
University of Alaska Anchorage\\
3211 Providence\\
Anchorage, AK 99508\\
smiley@math.uaa.alaska.edu}

\end{titlepage} 
\maketitle
\vskip -10pt
\centerline{\it In memory of Leonard Carlitz\rm}

\vskip 10pt

\begin{abstract}
The exponential generating functions $\sum_{n\geq 1}n^{n-m}{{z^n}\over{n!}}$ (for $m\in\mathbb{Z}$) are proven to be rational functions of the tree function ($\sum_{n\geq 1}n^{n-1}{{z^n}\over{n!}}$) and then of the endofunction function ($\sum_{n\geq 1}n^{n}{{z^n}\over{n!}}$), and some coefficients of these rational functions are identified with certain generalized Stirling, Eulerian, and Bernoulli numbers. Of the four families we cover, one had been previously analyzed by L. Carlitz.
\end{abstract}

\section{Introduction.}
In \cite{C}, Carlitz expressed $\sum_{n\geq 1}n^{n-m}{{z^n}\over{n!}}$ for $1\geq m\in\mathbb{Z}$ as a rational function of $\la (=T(z))=\sum_{n\geq 1}n^{n-1}{{z^n}\over{n!}}$. The coefficients in this rational function (written in lowest terms with factored parts) are positive integers which had previously appeared in \cite{M} and were subsequently given combinatorial interpretations in \cite{R1} and \cite{GS}. These numbers are now known as ``second-order Eulerian numbers'' (cf. \cite{GKP}).

In this note, we consider $m\geq 1$ and find polynomials in $\la$. We show that the coefficients of these polynomials are the virtual Stirling Numbers of the First Kind studied by D. Loeb in \cite{L1}. Previous tabulations and interpretations had appeared in \cite{DB} and elsewhere in the ``finite differences'' literature.

We introduce the variable $\zeta =\sum_{n\geq 1}n^n{{z^n}\over{n!}}={{\la}\over{1-\la}}$ and investigate the coefficients in the rational functions of $\zeta$ equivalent to those in the above sequence. In the case $m\leq 0$ we find polynomials of degree $2m+1$ which, in factored form, display integer coefficients equal to the ``associated Stirling Numbers of the second kind'' studied by Riordan in \cite{R2}. These have also been named ``second-order Stirling Numbers'' (in \cite{F}, for instance) and we show that easy manipulations produce new identities relating them to the second-order Eulerian numbers.

It remains to address the $\zeta$ expressions in the case $m\geq 1$. The rational functions, in lowest terms, are quotients of two polynomials each of degree $|m|$ in $\zeta$. While we obtain satisfactory results about the asymptotic behavior of the coefficients, the question of combinatorial interpretation is open.

\section{An Integral Operator.}

 For integer 
values of $m$ we define the sequence
$G_{m}(\la ):=R_m(z):=\sum_{n\geq 1}n^{n-m}{{z^n}\over{n!}}$. Thus, as defined above, $\la =G_1(\la )$ and $\zeta =G_0(\la )$. We define the formal integral operator $\mathbb{I}$ by $\mathbb{I}(F(\la ))=\int_0^\la F_m(\rho ){{1-\rho}\over{\rho}}\ d\rho$.

\begin{theorem}  $\mathbb{I}(G_m(\la ))=G_{m+1}(\la )$.
\end{theorem}
\vskip 0.1 true in
\noindent
\textit{Proof.} 
\begin{displaymath}G_{m+1}(\la )=R_{m+1}(z)=\int_0^zR_m(u){{du}\over{u}}=\int_0^\la R_m(\rho e^{-\rho}){{1-\rho}\over{\rho}}\ d\rho=\int_0^\la G_m(\rho ){{1-\rho}\over{\rho}}\ d\rho\qed
\end{displaymath}
\begin{corollary}  $G_{m-1}(\la )={{\la}\over{1-\la}}G_m'(\la )$.
\end{corollary}
\vskip 0.1 true in
\begin{comment}  The differential operator in the Corollary is essentially that used by Carlitz \cite{C}. Specifically he derived the following 
transform for nonnegative $m$: If
\begin{equation}
G_{-m}(\la ):=\sum_{n\geq 1}n^{n+m}{{z^n}\over{n!}}={{\la}\over{(1-\la )^{2m+1}}}\sum_{k=0}^m\left\langle\left\langle{{m}\atop{k}}\right\rangle\right\rangle\la^{k}
\end{equation}
\noindent
then
\begin{equation}
\left\langle\left\langle{{m}\atop{k}}\right\rangle\right\rangle = (k+1)\cdot\left\langle\left\langle{{m-1}\atop{k}}\right\rangle\right\rangle + (2m-k-1)\cdot\left\langle\left\langle{{m-1}\atop{k-1}}\right\rangle\right\rangle.
\end{equation}

Here the symbol $\l2e{{m}\atop{k}}\r2e$ has become known as the second order Eulerian number \cite[Table 270]{GKP}.  Combinatorial interpretations of  $\left\langle\left\langle{{n}\atop{k}}\right\rangle\right\rangle$ were subsequently given by Riordan\cite{R1} and Gessel and Stanley \cite{GS}.
\end {comment}
\begin{corollary} If $m>0$, $G_m(\la )$ is a polynomial of degree $m$ without constant term (Appendix 1). If we write $G_m(\la )=\sum_{k=1}^m(-1)^{k-1}g(m,k)\la^k$, then, for $m>1$, 
\begin{equation}
k\cdot g(m,k)=g(m-1,k)+g(m-1,k-1).
\end{equation} 
\end{corollary}
\vskip 0.1 true in
\noindent
We note that $g(m,1)=1$ and set $g(m,k)=0$ for $k>m>0$. These two boundary conditions, together with $(3)$, determine $g(m,k)$ for $m,k\geq 1$. A combinatorial interpretation of $g(m,k)$ is as follows: let $N=m-k+1$ and consider the set of all $k\times N$ matrices such that each column is a permutation of $\{1,2,\dots ,k\}$. Then $g(m,k)$ is the probability that, in a randomly selected matrix from this set, no one number is below the number $1$ in every column. 
\vskip 0.1 true in
\noindent
We note two previous appearances of the numbers $g(m,k)$ in the literature. In \cite[Appendix]{DB}, a definition of generalized Bernoulli numbers of negative degree is given in accordance with which $g(m,k)=-{{m-1}\choose{k}}B_{-m}^{(-k)}$. In \cite{L1}, the Stirling numbers of the First Kind are generalized such that $g(m,k)=(-1)^ks(-k,m-k)$. Many formulae for $g(m,k)$ may be found in these references and elsewhere (cf. Appendix 3).
 
\section{Pullback by Fractional Linear Transform $(m<0)$}
\vskip 0.3 true in
\noindent
If $\la$ may be called the tree variable, then $$\zeta =Z(z)=\sum_{n\geq 1}n^n{{z^n}\over{n!}}={{\la}\over{1-\la}}$$ is the endofunction variable. When we express
$R_m(z)=G_m(\la )$ in terms of $\zeta$ a beautiful phenomenon occurs: Carlitz' rational functions ($m$ nonpositive) become polynomials in $\zeta$ with integer coefficients strongly related to numbers called ``associated Stirling numbers of the second kind'' by Riordan [R2] or ``second-order Stirling numbers'' (denoted $\{\{{{m}\atop{k}}\}\}$) by Fekete[F]. We adopt the latter notation, in terms of which  $\left\{\left\{{{m}\atop{k}}\right\}\right\} =k\cdot\{\{{{m-1}\atop{k}}\}\} +(m-1)\cdot \{\{{{m-2}\atop{k-1}}\}\}$ is the defining 3-term recurrence relation.
\vskip 0.2 true in
\begin{propo} For $m\in\mathbb{N}$, 
\begin{equation}
G_{-m}(\la )=(1+\zeta )^{m+1}\sum_{k=1}^m\left\{\left\{ {{m+k}\atop{k}}\right\}\right\} \zeta^k:=H_{-m}(\zeta ).
\end{equation}  
\end{propo}
\vskip 0.1 true in
\noindent
\textit{Proof.} Making the substitution $\la ={{\zeta}\over{1+\zeta}}$ in $(1)$, and writing the Carlitz differential operator as $\zeta (1+\zeta )^2{{d}\over{d\zeta}}$, we calculate, with $G_{-m}(\la )=(1+\zeta )^{m+1}\sum_{k=1}^ma(m,k) \zeta^k$ and Theorem 1, that
$a(m,k)=(m+k-1)\cdot a(m-1,k-1)+k\cdot a(m-1,k)$, from which the result follows.$\qed$
\vskip 0.1 true in
These formal calculations have some pleasant concrete consequences.
\begin{corollary}
\begin{equation}
\sum_{k=0}^m\Bl2e{{m}\atop{k}}\Br2e (1+\zeta )^{m-k-1}\zeta^k=\sum_{k=1}^m\left\{\left\{ {{m+k}\atop{k}}\right\}\right\} \zeta^{k-1}
\end{equation}
\end{corollary}
Summation identities are produced by specializing $\zeta$. A pure counting identity follows from the binomial theorem:

\begin{corollary}
\begin{equation}
\left\{\left\{ {{n+q}\atop{q}}\right\}\right\} =\sum_{i=0}^n{{n-i-1}\choose{q-i-1}}\Bl2e{{n}\atop{i}}\Br2e.
\end{equation}
\end{corollary}
This is completed to an inverse pair of identities (as in \cite{R3}) by
\begin{corollary}
\begin{equation}
\Bl2e{{n}\atop{q}}\Br2e =\sum_{i=0}^n(-1)^{q-i}{{n-i-1}\choose{q-i}}\left\{\left\{ {{n+i+1}\atop{i+1}}\right\}\right\}
\end{equation}
\end{corollary}

\section{Pullback by Fractional Linear Transform $(m>0)$}
\vskip 0.3 true in
\noindent

The 
$\zeta -$expressions for positive $m$ may be written $H_m(\zeta )=(1+\zeta )^{-m}\sum_{k=1}^mh(m,k)\zeta^k$. The rational numbers $h(m,k)$ do not appear, to our knowledge, in the literature (see Appendix 2, top). In particular, combinatorial interpretation of $h(m,k)$ is unknown. It is possible, of course, to use our operators and the correlative properties of $g(m,k)$ to analyze these coefficients quantitatively. We state the basic facts without proof.
\vskip 0.1 true in
\noindent
First we note that these coefficients may be given by the formula
\begin{equation}
h(m,k)=\sum_{j=1}^k(-1)^{j-1}{{m-j}\choose{k-j}}g(m,j)
\end{equation}
\noindent
which may be combined with known formulas for $g(m,k)$. A better clue toward finding a combinatorial interpretation may come from the rather simple recurrence relation
\begin{equation}
k \cdot h(m,k) = h(m-1,k) + (m-k+1)\cdot h(m,k-1).
\end{equation}
\noindent
Since $g(m,k)$ has an easy interpretation as a probability, one might hope to interpret $h(m,k)$ as an expected value in some finite problem. To this end we note that both $1!2!\cdots m!\cdot h(m,k)$ and $(1!2!\cdots k!)^{m-k+1}\cdot h(m,k)$ are integral. We further state
\begin{equation}
\sum_{j=1}^m(-1)^{j-1}h(m,j)=1/m!
\end{equation}
\noindent
These facts may be easily proved by induction. It is also easy to verify that 
\begin{equation}
\lim_{m\to\infty}h(m,m)=1/e
\end{equation}
\vskip 0.1 true in
The growth of $h(m,k)$ as $m\to\infty$ for fixed $k$ is of interest. Experiment strongly suggests that $(k-1)!\cdot h(m,k)$ is asymptotic to a polynomial $p_k(m)$ of degree $k-1$ in $m$ with nonnegative integer coefficients.

\section{Conclusion.}
\vskip 0.3 true in
\noindent
The variables $\la =T(z)$ and $\zeta$, as well as $W:=-T(-z)$ have been used 
extensively in analyses of random mappings and random graphs (se \cite{JKLP,FO} for surveys) in which other $R_m(z)$ appear. The latter also occur naturally as exponential generating functions of tree-like structures. We hope that our very basic results may be useful to further such studies.
\newpage
\centerline{\textbf{Appendix 1: The $\la$ Sequence}}
 For $\la (=T(z))=\sum_{n\geq 1}n^{n-1}{{z^n}\over{n!}}$ and integer 
values of $m$ we define the sequence
$G_{m}(\la ):=R_m(z):=\sum_{n\geq 1}n^{n-m}{{z^n}\over{n!}}$. Then

$$\dots$$
$$-{{1}\over{720}}\la^6+{{137}\over{7200}}\la^5-{{415}\over{3456}}\la^4+{{575}\over{1296}}\la^3-{{31}\over{32}}\la^2+\la\leqno{G_6(\la )=}$$
$${{1}\over{120}}\la^5-{{25}\over{288}}\la^4+{{85}\over{216}}\la^3-{{15}\over{16}}\la^2+\la\leqno{G_5(\la )=}$$
$$-{{1}\over{24}}\la^4+{{11}\over{36}}\la^3-{{7}\over{8}}\la^2+\la\leqno{G_4(\la )=}$$
$${{1}\over{6}}\la^3-{{3}\over{4}}\la^2+\la\leqno{G_3(\la )=}$$
$$-{{1}\over{2}}\la^2+\la\leqno{G_2(\la )=}$$
$$\la\leqno{G_1(\la )=}$$
$${{\la}\over{1-\la}}\leqno{G_0(\la )=}$$
$${{\la}\over{(1-\la )^3}}\leqno{G_{-1}(\la )=}$$
$${{\la (1+2\la )}\over{(1-\la )^5}}\leqno{G_{-2}(\la )=}$$
$${{\la (1+8\la +6\la^2)}\over{(1-\la )^7}}\leqno{G_{-3}(\la )=}$$
$${{\la (1+22\la +58\la^2+24\la^3)}\over{(1-\la )^9}}\leqno{G_{-4}(\la )=}$$
$${{\la (1+52\la +328\la^2+444\la^3+120\la^4)}\over{(1-\la )^{11}}}\leqno{G_{-5}(\la )=}$$
$$\dots$$
\eject
\centerline{\textbf{Appendix 2: The $\zeta$ Sequence}}
 For  $\zeta =\sum_{n\geq 1}n^n{{z^n}\over{n!}}$ and integer 
values of $m$ we define the sequence
$H_{m}(\zeta ):=R_m(z):=\sum_{n\geq 1}n^{n-m}{{z^n}\over{n!}}$. Then

$$\dots$$
$${{{\zeta}(1+{{129}\over{32}}{\zeta}+{{8513}\over{1296}}{\zeta}^2+{{691}\over{128}}{\zeta}^3+{{96547}\over{43200}}{\zeta}^4+{{96547}\over{259200}}{\zeta}^5)}\over{(1+{\zeta})^6}}\leqno{H_6(\zeta )=}$$
$${{{\zeta}(1+{{49}\over{16}}{\zeta}+{{1547}\over{432}}{\zeta}^2+{{1631}\over{864}}{\zeta}^3+{{1631}\over{4320}}{\zeta}^4)}\over{(1+{\zeta})^5}}\leqno{H_5(\zeta )=}$$
$${{{\zeta}(1+{{17}\over{8}}{\zeta}+{{14}\over{9}}{\zeta}^2+{{7}\over{18}}{\zeta}^3)}\over{(1+{\zeta})^4}}\leqno{H_4(\zeta )=}$$
$${{{\zeta}(1+{{5}\over{4}}{\zeta}+{{5}\over{12}}{\zeta}^2)}\over{(1+{\zeta})^3}}\leqno{H_3(\zeta )=}$$
$${{{\zeta}(1+{{1}\over{2}}{\zeta})}\over{(1+{\zeta})^2}}\leqno{H_2(\zeta )=}$$
$${{\zeta}\over{1+\zeta}}\leqno{H_1(\zeta )=}$$
$${\zeta}\leqno{H_0(\zeta )=}$$
$$\zeta (1+\zeta )^2\leqno{H_{-1}(\zeta )=}$$
$$\zeta (1+\zeta )^3(1+3\zeta )\leqno{H_{-2}(\zeta )=}$$
$$\zeta (1+\zeta )^4(1+10\zeta +15\zeta^2)\leqno{H_{-3}(\zeta )=}$$
$$\zeta (1+\zeta )^5(1+25\zeta +105\zeta^2+105\zeta^3)\leqno{H_{-4}(\zeta )=}$$
$$\zeta (1+\zeta )^6(1+56\zeta +490\zeta^2+1260\zeta^3+945\zeta^4)\leqno{H_{-6}(\zeta )=}$$

$$\dots$$
\eject
\vglue -1.2 true in
\centerline{\textbf{Appendix 3:  Extended Remark on $g(m,k)$}}
\vskip 0.3 true in
The tabulation of  $\Bl2e{{m}\atop{k}}\Br2e$ led Carlitz to ask for combinatorial interpretations. In contrast, the coefficients $g(m,k)$ are well known. In \cite{DB}, a definition of generalized Bernoulli numbers of negative degree is given from which we see that $g(m,k)=-{{m-1}\choose{k}}B_{-m}^{(-k)}$. Also $g(m,k)=(-1)^ks(-k,m-k)$, a virtual Stirling number as defined in \cite{L1}.
\par
\vskip 0.3 true in
\noindent
Consider the recurrence relation
$$k\cdot g(m,k)=g(m-1,k)+g(m-1,k-1)$$
\noindent
for $m,k\geq 1$ where we initialize by $g(m,k)=0$ if $k>m$, and $g(m,1)=1$. Since this is our definition, the verification of the following formulas is routine.
\noindent
The unique solution may be expressed as a sum of Egyptian Fractions:
$$g(m,k)={{1}\over{k!}}\sum_{1\leq i_1\leq i_2\leq\dots\leq i_{m-k}\leq k}{{1}\over{i_1i_2\cdots i_{m-k}}}$$
\noindent
An alternative form is related to an iterated difference of $f(x)={{1}\over{x^{m-k+1}}}$:
$$g(m,k)={{1}\over{k-1!}}\sum_{j=0}^{k-1}{{k-1}\choose{j}}{{(-1)^j}\over{(j+1)^{m-k+1}}}$$
\noindent
From the recurrence for the $g(m,k)$'s also follow generating functions:
$$g(m,k)=[z^{m-k}]\ \prod_{p=1}^k{{1}\over{p-z}}$$
\noindent
which may be written in a connection coefficient formula:
$${{1}\over{(x+1)(x+2)\cdots (x+k)}}=\sum_{n=0}^\infty (-1)^ng(k+n,k)x^n$$
\noindent
(linking $g(m,k)$ with the signless Stirling cycle numbers).
\vskip 0.3 true in
An integral of a polynomial over the hypercube yields the same numbers:
$$g(m,k)={{1}\over{(k-1)!}}\int_0^1\int_0^1\cdots\int_0^1(1-x_1x_2\cdots x_{m-k+1})^{k-1}dx_1\ dx_2\cdots dx_{m-k+1}$$
\noindent
since it evaluates to the iterated difference form above after applying the binomial theorem. 
\vskip 0.3 true in
The rational number $g(m,k)$ has the natural denominator $(k!)^{m-k+1}$. The corresponding numerators $N(m,k)$ form a triangle of integers generated by
$$N(m,k)=[z^{m-k}]\ \prod_{p=1}^k{{1}\over{1-{{k!}\over{p}}z}}$$
\noindent
and previously tabulated and studied in  David, Kendall, and Barton, \textit{Symmetric Function and Allied Tables}, Cambridge Univ. Press, London, 1966 (Table 5.4.1).

\eject

\vskip 0.4 true in
\leftline{\textbf{Web Resource}}
\vskip 0.3 true in
Some additional data on the questions raised in Section 4 may be found at
\par
\noindent
\texttt{http://www.saturn.math.uaa.alaska.edu/$\tilde{\ }$smiley/numbertriangle.html}

\end{document}